\documentclass[10pt]{article}
\usepackage[latin1]{inputenc}
\usepackage{amsmath}
\usepackage{amsfonts}
\usepackage{amssymb}
\newcommand{\ZA}{{\mathcal A}}

\newcommand{\ZEP}{\epsilon}

\newcommand{\zg}{\gamma}

\newcommand{\intT}{\int_0^T}
\newcommand{\intt}{\int_0^t}
\newcommand{\ints}{\int_0^s}
\newcommand{\intr}{\int_0^r}

\newtheorem{Theorem}{Theorem}

\newtheorem{Remark}[Theorem]{Remark}

\newtheorem{Definition}[Theorem]{Definition}

\newcommand{\zaa}{\alpha}

\newcommand{\ZDE}{\delta}
\newcommand{\zt}{\tau}
\newcommand{\zdia}{~~\rule{1mm}{2mm}\par\medskip}

\newcommand{\ZLA}{\label}
\newcommand{\ZIN}{\infty}
\newcommand{\zProof}{{\bf\underbar{Proof}.}\ }

\newcommand{\ZBI}{\bibitem}
\newcommand{\ZD}{\;\mbox{\rm d}}
\newcommand{\zl}{\lambda}
\newcommand{\ZSI}{\sigma}
\newcommand{\ZOMq}{\Omega}
\newcommand{\ZSUno}{\sum _{n=1}^{+\ZIN}}
\newcommand{\Mo}{\mathbb{M}}

\newcommand{\bbZ}{\mathbb{Z}}

\author{L. Pandolfi}
\title{Cosine operator and controllability of the wave equation with memory revisited}
\begin{document}
\maketitle

\section{Introduction}

In this paper we consider the following equation with memory:
\begin{equation}
\ZLA{eq:sistemaVERS4}
w'=2\zaa w(t)+\intt N(t-s)\Delta w(s)\ZD s
\end{equation}
where $w=w(x,t)$ with $t>0$ and  $x\in\ZOMq$, a region with $C^2$ boundary and $N\in H^3(0,T)$ for every $T>0$.
We associate the initial and boundary conditions
\[
\left\{\begin{array}
{l}
w( 0)=w_0 \in L^2(\ZOMq)\,,\\
 w( t)=f( t)\ {\rm if}\ x\in\Gamma\subseteq\ZOMq\,,\quad w( t)=0\ {\rm if}\ x\in\partial\ZOMq\setminus\Gamma
\end{array}
\right.
\]
($\Gamma$ is relatively open in $\partial\ZOMq$. The case $\Gamma=\partial\ZOMq$ is not excluded).

The function $f$ is a control which we use to steer the initial datum $w_0\in L^2(\ZOMq)$ to a target $\xi\in L^2(\ZOMq)$ at a certain time $T$. This kind of control problem has been studied by several authors and with different methods, since a system of the form~(\ref{eq:sistemaVERS4}) is important for the applications in viscoelasticity, thermodynamics of materials with memory and nonfickian diffusion. Note that in viscoelasticity also controllability of the pair $(w,w')$ of the deformation and velocity has to be studied but here for simplicity we confine ourselves to the controllability of the sole component $w $. It is a fact that: 1) the controllability of the sole component $w$ is sufficient for the solutions of source identification problems,
 see~\cite{PandDCDS2,PandViscoUnderTRACT,PANDsurvey}; 2) also the controllability of the pair of the deformation and the stress (or the flux) has its interest, and this is a new problem which appears in the case of systems with memory, see~\cite{AvdoninPANDOLFIflux1,AvdoninPANDOLFIflux2,PandolfiTRACTION-deformation}.

The key idea  which underlines essentially all the papers on controllability of Eq.~(\ref{eq:sistemaVERS4}) is that the controllability properties of the \emph{associated wave equation}
\begin{equation}
\ZLA{eq:associatedWAVEvers4}
u''=\Delta u+F\,,\qquad \left\{\begin{array}{l}
u(0)=u_0\,,\ u'(0)=u_1\,,\\
u( t)=f( t)\ {\rm if}\ x\in\Gamma\subseteq\ZOMq\,,\quad u( t)=0\ {\rm if}\ x\in\partial\ZOMq\setminus\Gamma
\end{array}\right.
\end{equation}
can be lifted to the system with memory~(\ref{eq:sistemaVERS4}).

The paper~\cite{PandAMO} proved that the control properties of~(\ref{eq:associatedWAVEvers4}) can be lifted to the system~(\ref{eq:sistemaVERS4}) using cosine operator theory.
Here  we intend to revise and improve this approach.

The organization of the paper is as follows: first we combine the MacCamy trick to give a definition/representation of the solutions of 
Eq.~(\ref{eq:sistemaVERS4}), see section~\ref{Sect:SOLUmemory-VERS-4}. 
In particular we prove that for every $f\in L^2(0,T;L^2(\Gamma))$ and every $w_0\in L^2(\ZOMq)$ a solution exists, such that $w\in C([0,T];L^2(\ZOMq))$ for every $T>0$. This justify the following definition of the reachable set at time $T$ (the index $M$ is for ``memory''):
\[
R_M(T)=\{ w(T)\,,\quad f\in L^2(0,T;L^2(\Gamma))\}\,.
\]  

The final result is the proof that if the wave equation is controllable then system~(\ref{eq:sistemaVERS4}) is controllable too (see the precise statement in Theorem~\ref{teo:CONTROLLprovatoDAVVERO}). This is in two steps: in the first step we prove that $R^\perp (T)$ is finite dimensional and then we prove that its orthogonal is  reduced to the subspace $0$.

These arguments depend on known properties of the wave equation, which are recalled in Section~\ref{SECTpropreWAVE-4}.

\subsection{\ZLA{sec:commentPREVIOUSreferences-4}Comments on previous results}

Controllability of Equation~(\ref{eq:sistemaVERS4}) has been studied using different methods which are reviewed in~\cite{PandLIBRO}. The papers~\cite{LeugeMOME} uses Fourier expansions and moment methods, an approach extended in~\cite{PandolfiSHARP} (see also~\cite{PandLIBRO} and references therein). Extension to~(\ref{eq:sistemaVERS4}) of the inverse inequality of the wave 
equation is in~\cite{Kim1993}
while Carlemn estimates are used in~\cite{FuZHANG1}. Here we extend and improve the operator approach in~\cite{PandAMO}.

\section{\ZLA{SECTpropreWAVE-4}\ZLA{Sect:proprieWAVE}The properties of the wave equation}

%\section{\ZLA{Sect:proprieWAVE}Properties of the wave equation}
We need few pieces of  information on the wave equation~(\ref{eq:associatedWAVEvers4}). 
We introduce the operators $A$, $\ZA$ and $D$:
 
\[ 
\ZA=i(-A)^{1/2}\quad \mbox{where \  $ {\rm dom}\, A=H^2(\ZOMq)\cap H^1_0 (\ZOMq)) $, $ A\phi=\Delta\phi $}
 \]
 while   the operator $ D $, the Dirichlet operator,   is defined by
  \[ 
  u=Df\ \iff\ \Delta u=0\, \quad  
u(x)= f(x)  \ {\rm on}\ \Gamma\,,\ u=0\ {\rm on}\ \partial\ZOMq\setminus\Gamma\,.
   \]
The operator $\ZA$ generates a strongly continuous group, so that we can consider the strongly continuous  operators $R_+(t)$ and $R_-(t)$ defined by:
\[
R_+(t)=\frac{1}{2}\left [ e^{\ZA t}+e^{-\ZA t}\right ]\,,\qquad R_-(t)=\frac{1}{2}\left [ e^{\ZA t}-e^{-\ZA t}\right ]\,.
\]
The operator $R_+(t)$ is the \emph{cosine operator} generated by $A$ and its key property is
\[
R_+(t)R_+(\zt)=\frac{1}{2}\left [R_+(T+\zt)+R_+(t-\zt)\right ]\,.
\]
This equality holds for every real  $t$ and $\zt$.

Let $u_0\in L^2(\ZOMq)$, $u_1\in H^{-1} (\ZOMq)$, $F\in L^1(0,T;L^2(\ZOMq))$ and $f\in L^2(0,T;L^2(\Gamma))$.
It is known that problem~(\ref{eq:associatedWAVEvers4}) admits a unique solution $u\in C([0,T];L^2(\ZOMq))\cap C^1([0,T];H^{-1} (\ZOMq))$ which is given by
\begin{equation}\ZLA{eq:SOLUequaONDE}
\begin{array}{l}
\displaystyle u(t)=R_+(t) u_0+\ZA^{-1}R_-(t)u_1+\ZA^{-1}\intt R_-(t-s)F(s)\ZD s-\\
\displaystyle -\ZA\intt R_-(t-s)Df(s)\ZD s\,,\\[2mm]
\displaystyle
u'(t)= \ZA R_-(t)u_0+R_+(t) u_1+\intt R_+(t-s)F(s)\ZD s-A\intt R_+(t-s)D f(s)\ZD s\,.
\end{array}
\end{equation}

The following  result is known (see~\cite{LIONSlibro}). Let $\zg_1$ be the exterior normal derivative,
\[
\zg_1\phi(x)=\frac{\partial}{\partial n}\phi(x)\,,\quad x\in \partial\ZOMq\,.
\]

\begin{Theorem}
The following properties hold for the \emph{memoryless} wave equation~(\ref{eq:associatedWAVEvers4}). We state separately the effects of $ u_0 $, $ u_1 $, $ F $ and    of the boundary control $f$.
\begin{enumerate}
\item
Let $f=0$ and $u_0\in H^1_0 (\ZOMq)$, $u_1\in L^2(\ZOMq)$, $F\in L^1(0,T;L^2(\ZOMq))$. Then $u(t)\in C([0,T];H^1_0 (\ZOMq))\cap C^1([0,T];L^2(\ZOMq))$ and it is a linear and continuous function of $u_0$, $u_1$, $F$    in the specified spaces.  Furthermore,
  for every $T>0$ there exists $M>0$ such that
\begin{equation}\ZLA{eq:DIREdisegPUREonde}
\int_\Gamma \int_0^T\left  |\zg_1 u(t) \right |^2\ZD t\, \ZD \Gamma \leq M\left [ \|u_0\|^2 _{H^1_0 (\ZOMq)}+\|u_1\| ^2_{L^2(\ZOMq)}+\|F\|^2 _{L^1(0,T;L^2(\ZOMq))}\right ]\,.
\end{equation}

\item
If $f=0$ and $u_0\in L^2(\ZOMq)$, $u_1\in H^{-1} (\ZOMq)$, $F\in L^1(0,T;L^2(\ZOMq))$ then $u(t)\in C([0,T];L^2(\ZOMq))\cap C^1([0,T];H^{-1} (\ZOMq))$ and it is a linear and continuous function of $u_0$, $u_1$, $F$    in the specified spaces. 

\item if $f\in L^2(0,T;L^2(\Gamma)) $ and $u_0=0$, $u_1=0$, $F=0$ then $u(t)\in C([0,T];L^2(\ZOMq))\cap C^1([0,T];H^{-1} (\ZOMq))$ and depends continuously on $f$.
\end{enumerate}
\end{Theorem}
 
The previous properties justify the following definition, where the control time is called $ 2T $ for later convenience:
\begin{Definition}\ZLA{eq:defiCONTROwave}
The wave equation~(\ref{eq:associatedWAVEvers4}) is controllable at time $2T$ if for every $u_0$ and $\xi$ in $L^2(\ZOMq)$, $u_1$ and $\eta$ in $H^{-1} (\ZOMq)$ and $F\in L^1(0,T;L^2(\ZOMq))$ there exists a control $f\in L^2(0,T;L^2(\Gamma))$ such that
\[
u(2T)=\xi\,,\quad u'(2T)=\eta\,.
\]
\end{Definition}
It is known that
\begin{enumerate}
\item if $\Gamma$ is ``too small'' then there exists no time   at which the wave equation is  controllable.
\item there exist subset $\Gamma$ of $\partial\ZOMq$ (for example, $\Gamma=\partial\ZOMq$) such that controllability holds for a suitable time.
\item if controllability holds at time $2T$ then it holds also at every larger time.
\item controllability does not depend on $u_0$, $u_1$ and $F$ so that when studying controllability we can assume $u_0=u_1=0$, $F=0$. So, controllability is the property that the following map is surjective. The map acts from $ L^2(0,T;L^2(\Gamma)) $ to $ L^2(\ZOMq)\times H^{-1}(\ZOMq) $ and it is defined by
\begin{eqnarray*}
&&f\mapsto \Lambda_0(2T) f=\left ( \Lambda_0^1(2T), \Lambda_0^2(2T)\right )f=\\
&&=
 \left (
\ZA\int_0^{2T} R_-(2T-s)Df(s)\ZD s\,,\ A\int_0^{2T} R_+(2T-s)D f(s)\ZD s
 \right )\,.
\end{eqnarray*}
\end{enumerate}

\subsection{A consequence in terms of bases}

  It is known that the operator $A$ is selfadjoint with compact resolvent. Hence, $L^2(\ZOMq)$ has an orthonormal basis whose elements $\phi_n(x)$ are eigenvectors of $A$:
  \[
A\phi_n= -\zl_n^2\phi_n\,.  
  \]
 It is a fact that $ \zl_n^2>0$ hence $\zl_n$ is real and we can choose $\zl_n>0$. The eigenvalues are not distinct, but the eigenvectors with the same eigenvalue are finite in number.
 
 The operators $ R_+(t) $ and $ R_-(t) $ have a simple representation in terms of $ \phi_n(x) $:
  
 \begin{eqnarray*}
 && R_+(t)\left (\ZSUno c_n\phi_n(x)\right )=\ZSUno \phi_n(x)\left ( c_n\cos\zl_n t\right )\,,\\
  && R_-(t)\left (\ZSUno c_n\phi_n(x)\right )=\ZSUno \phi_n(x)\left ( c_n\sin \zl_n t\right )\,. 
  \end{eqnarray*}
 Furthermore, we know that  (see~\cite[Prop.~10.6.1]{TuksnakWeiss} and note that our operator $ A  $ is   $ -A_0$    in~\cite{TuksnakWeiss})
 \[ 
 \int_\ZOMq \phi_n(x)Df\ZD x=-\frac{1}{\zl_n^2}\int_\Gamma (\zg_1\phi_n) f\ZD \Gamma\,.
  \]
So, $ -\Lambda_0(2T) f $ has the following ``concrete'' representation:
\[
 \begin{array}{l}
\displaystyle 
\left (
 \ZSUno\phi_n(x) 
\int_0^{2T}\int_\Gamma\left ( \frac{\zg_1\phi_n}{\zl_n}\right )\left ( \sin\zl_n s\right ) f(x,2T-s)\ZD\Gamma\,\ZD s\,,
 \right.
 \\
 \displaystyle
 \left.   
 \ZSUno\left (\zl_n \phi_n(x)\right )
 \int_0^{2T} \int_\Gamma \left (\frac{\zg_1\phi_n}{\zl_n}\right )\left ( \cos\zl_n s\right ) f(x,2T-s)\ZD\Gamma\,\ZD s
\right )\,.
 \end{array}
 \]
 
 It is known that $ \{\zl_n\phi_n\}  $ is an orthonormal basis of $ H^{-1} (\ZOMq) $ (in an inner product whose norm is equivalent to the standard norm) so that 
  every target has the representation
 \[ 
\xi=\ZSUno \xi_n \phi_n\,,\quad   \eta = \ZSUno  \eta_n\left (\zl_n\phi_n\right )   \quad \{\xi_n\}\in l^2\,,\ \{\eta_n\}\in l^2\,.
  \]
 Controllability  is equivalent to the solvability of the following moment problem, in terms of a \emph{real} function $ f $:
  \begin{equation}
\ZLA{eq:MOMRondeREALE}
\int_0^{2T}\int_\Gamma \Psi_n e^{i\zl_n s}f(x,2T-s)\ZD \Gamma\,\ZD s=\eta_n+i\xi_n=c_n\,,\quad n\in\mathbb{N}\qquad \Psi_n=\frac{\zg_1\phi_n}{\zl_n}\,.
\end{equation}
Note that $ \{c_n\} $ is an \emph{arbitrary complex valued $ l^2 $-sequence.} 

We introduce  
\[ 
\bbZ'=\bbZ\setminus\{0\}\,,\quad  \zl_n=-\zl _{-n}\,,\ \phi_n=\phi _{-n} \ {\rm for}\ n<0\;.
 \]
Then, the   moment problem~(\ref{eq:MOMRondeREALE}) is equivalent to
   \begin{equation}
\ZLA{eq:MOMRondeCOMPLE}
\int_0^{2T}\int_\Gamma \Psi_n e^{i\zl_n s} h(x,2T-s)\ZD \Gamma\,\ZD s =c_n\,,\quad n\in\bbZ' 
\end{equation}
where now $ h\in L^2(0,2T;L^2(\Gamma)) $ is \emph{complex valued } and $\{c_n\}\in l^2(\bbZ')$ is arbitrary.

We introduce the moment operator    
\[ 
\Mo_0 h=\int_0^{2T}\int_\Gamma \Psi_n e^{i\zl_n s} h(x,2T-s)\ZD \Gamma\,\ZD s\,.
 \] 
 The fact that the transformation $ f\mapsto \Lambda_0(2T) f $ is continuous from $ L^2(0,2T;L^2(\Gamma)) $ to $ L^2(\ZOMq)\times H^{-1} (\ZOMq) $ implies that $ \Mo_0\in \mathcal{L}\left (L^2(0,2T;L^2(\Gamma)), l^2(\bbZ')\right ) $.
  So, its domain is   $ L^2(0,2T;L^2(\Gamma)) $  and its restiction to $ {\rm cl\, span} \{\Psi_ne^{i\zl_n t}\} $ is invertible.
 
 If the wave equation is controllable at time $ 2T $ then $ \Mo_0 $ is \emph{surjective,} so that its inverse (as an operator from $ {\rm cl\, span} \{\Psi_ne^{i\zl_n t}\} \subseteq L^2(0,2T;L^2(\Gamma))$ to $ l^2(\bbZ') $) is bounded. This implies (see~\cite[p.~22]{LoretiKOMORNIK} and~\cite{PandolfiSHARP,PandLIBRO}):
 \begin{Theorem}\ZLA{teo:controMETAonde}
 Let the associated wave equation be controllable in time $2T$.
 Then:
\begin{itemize}
 \item  
  the sequence $  \{\Psi_ne^{i\zl_n t} \}_{n\in \bbZ'}$ is a Riesz sequence in $L^2(0,2T;L^2(\Gamma))$;
  \item 
  the sequences $  \{\Psi_n\cos\zl_n t\}_{n\in \mathbb{N}}$, $  \{\Psi_n\sin\zl_n t\}_{n\in \mathbb{N}}$ are Riesz sequences in $ L^2(0,T;L^2(\Gamma)) $;
  \item  the operator $\Lambda_0^1(T)$ is surjective.

 \end{itemize}     
 \end{Theorem}

\section{\ZLA{Sect:SOLUmemory-VERS-4}The solutions of System~(\ref{eq:sistemaVERS4})}

Different methods have been proposed to study the solutions of Eq.~(\ref{eq:sistemaVERS4}). Here  we follow a method based on the use of cosine operator theory.

We apply a transformation,  known as MacCamy trick, to the solutions of Eq.~(\ref{eq:sistemaVERS4}). This is a formal step, since the solutions are not yet defined. Formally, computing the derivatives of both the sides of~(\ref{eq:sistemaVERS4}) we get
\begin{equation}
\ZLA{eq:sistemaVERS4-DERTIVATA}
w''=2\zaa w'=\left [ \Delta w+\intt N'(t-s)\Delta w(s)\ZD s\right ]+F(t)
\end{equation}
(the affine term should be zero. We inserted it here since even if $F=0$ an affine term will appear in the following computations).

We consider~(\ref{eq:sistemaVERS4-DERTIVATA})     as a Volterra integral equation in the unknown   $ \Delta w $. Let $ \tilde M(t) $ be the resolvent kernel of $ M(t) =N'(t)$, i.e. $\tilde M(t)$ is the unique solution of
\[
\tilde M(t)=M(t)-\intt N'(t-s)\tilde M(s) \ZD s\,.
\]
Then we get formally 
\[
\Delta w(t)=w''(t)-2\zaa w'-F -\intt  \tilde M(t-s)\left [ w''(s)-2\zaa w'(s)-F(s)\right ]\ZD s\,.
\]
In this equation, $w(0)=w_0$ and $w'(0)=0$. For the following it is convenient to study this equation in the more general case
\[
w'(0)=w_1\,,
\]
possibly different from zero.
We integrate by parts and we find a system of the following form:
\begin{equation}
\ZLA{eq:dopoMacCamyANTE}
w''(t)=\Delta w(t)+a w'(t)+b w(t)+\intt \tilde M_1(t-s) w(s)\ZD s+F_1(t)
\end{equation}
where $a$ and $b$ are suitable constants and
\[
F_1(t)=F(t)-\intt \tilde M(t-s)F(s)\ZD s-\tilde M(t) w_1- \tilde M'(t)w_0\,.
\]
Note the dependence of $F_1(t)$ on  the initial conditions $w_0$ and $w_1$ and note also that after the MacCamy trick the laplacian does not appear  in the memory of the integral.

It is simple to see that the transformation $w(t)\mapsto e^{-at/2} w(t)$ can be used to remove the velocity term from~(\ref{eq:dopoMacCamyANTE}) (this changes $\tilde M (t)$ to $M(t)=e^{-2at/2}\tilde M (t)=K(t)$ and similar transformation of $F_1(t)$ and the boundary control $f$). So, we study the problem
\begin{equation}
\ZLA{eq:dopoMacCamy}
\begin{array}{l}
w''(t)=\Delta w(t) +bw(t) +\intt   K (t-s) w(s)\ZD s+F (t)\,,\\[5pt]
% \left\{
% \begin{array}{l}
w(0)=w_0\,,\ w'(0)=w_1\\
w=f\ {\rm on}\ \Gamma\,,\ w=0\  {\rm on}\ \partial\ZOMq 
%\end{array}
%\right.
\end{array}
\end{equation} 
where $F_1(t)$ has been renamed  $F(t)$ and it is a continuous affine function of $w_0$ and $w_1$.

 This is a perturbed wave equation and we can use formula~(\ref{eq:SOLUequaONDE}) in order to get a Volterra integral equation for $w(t)$. It is convenient to write separately the formula for the contribution of the boundary control $f$ and for the contribution of $w_0$, $w_1$ and $F$:
  
\begin{eqnarray}
\nonumber
%\left\{\begin{array}{ll}
  w(t)&=&R_+(t) w_0+\ZA^{-1}R_-(t)w_1+\ZA^{-1}\intt R_-(t-s)F(s)\ZD s+\\
\ZLA{soluMEMOnoCONTROLLO}
 &&+\ZA^{-1}\intt R_-(t-s) \left [bw(s)+\ints K(s-r) w(r)\ZD r\right ]\ZD s\,,\\[3mm]
%w'(t)=&\ZA R_-(t)w_0+R_+(t) w_1+\intt R_+(t-s)F(s)\ZD s+\\
%\displaystyle&+\intt R_+(t-s)\left [bw(s)+\ints K(s-r) w(r)\ZD r\right ]      \ZD s
%\end{array}\right.
%\end{equation}
% When $w_0=0$, $w_1=0$, $F=0$ but $f\neq 0$ then we have:
%  
%\begin{equation}
\nonumber
%\left\{\begin{array}{ll}
%\displaystyle
w(t)&=&-\ZA\intt R_-(t-s)Df(s)\ZD s+\\
\ZLA{eq:diWconCONTROLLO}
&&+\ZA^{-1}\intt R_-(t-s)\left [bw(s)+\ints K(s-r)w(r)\ZD r\right ]\ZD s\,.
%\displaystyle
%w'(t)=&-A\intt R_+(t-s) Df(s)\ZD s+\\
%\displaystyle&+\ZA^{-1}\intt R_-(t-s)\left [ bw(s)+K(0) w(s)+\ints K'(s-r) w(r)\ZD r\right ]\ZD s\,. 
%\end{array}\right.
\end{eqnarray}
 The general solutions is the sum of the two but for our applications we keep distinct the two formulas.
 
 Both the formulas~(\ref{soluMEMOnoCONTROLLO}) and~(\ref{eq:diWconCONTROLLO}) for $ w(t) $ have the following general form
 \begin{equation}
\ZLA{eq:formaGENw}
w(t)= u(t)+ \ZA ^{-1}\intt L(t-s) w(s)   \ZD s\,,\quad L(t)w=b R_-(t)w+\intt K(t-r)R_-(r)w\ZD r\,.
\end{equation}

Note that $ u(t)  $ solves the associated wave equation.

 This is a Volterra integral equation for $ w(t) $ which we solve using Picard iteration:
 \begin{align}
\nonumber w(t)&= u(t)+\ZA^{-1}\intt L(t-s)w(s)\ZD s=\\
\ZLA{eq:RAPPREwCOMEconvoluzione}&= u(t)+\ZA^{-1}\intt L(t-s) u(s)\ZD s+\sum _{k=2}^{+\ZIN} \left ( \ZA^{-1}\right )^k L^{*k}*u 
\end{align}
where the exponent $ {~}^{*k} $ denotes iterated convolution.

We introduce the kernel $ H(t) $:
\[ 
H(t)=\sum _{k=1}^{+\ZIN} \left (\ZA^{-1}\right )^k\left (L\right )^{*k}=
\ZA^{-1}\left (\sum _{k=1}^{+\ZIN} \left (\ZA^{-1}\right )^{k-1}\left (L^{*k}\right ) \right ) 
 \]
 so that
\[ 
w(t)=u(t)+\intt H(t-s)u(s)\ZD s\;.
 \]
 When the boundary control is $ f=0 $ this formula specializes to
  \begin{equation}\ZLA{equA:FORMu;adiW}
 \begin{array}{ll}
\displaystyle w(t)=&R_+(t) w_0+\ZA^{-1}R_-(t)w_1+\ZA^{-1}\intt R_-(t-s)F(s)\ZD s+\\[2mm]
&+ \intt H(t-s) \left [R_+(s) w_0+\ZA^{-1}R_-(s)w_1+
\ZA^{-1}\ints R_-(s-r)F(r)\ZD r\right ]\ZD s
\end{array} 
\end{equation}
while the corresponding formula with $ w_0=w_1=0 $ and $ F=0 $ is
\begin{equation}
\ZLA{equA:FORMu;adiW-CONTROLLO}
 w(t)=-\ZA\intT R_-( t-s)Df(s)\ZD s-\intt  H(t-s)\ZA\int_0^s R_-(s-r) Df(r)\ZD r\,\ZD s
\end{equation}

The properties of the solutions of the wave equation that we recalled in Sect.~\ref{Sect:proprieWAVE} imply:
\begin{Theorem}
\ZLA{teO:proprietaEQUAconMEMORIA}
Let $F\in L^1(0,T;L^2(\ZOMq))$, $f\in L^2(0,T;L^2(\Gamma))$, $w_0\in L^2(\ZOMq)$ and $w_1\in H^{-1} (\ZOMq)$. Then $w\in C([0,T];L^2(\ZOMq))\cap C^1([0,T];H^{-1} (\ZOMq))$. If $f=0$ and $w_0\in H^1_0 (\ZOMq)$, $w_1\in L^2(\ZOMq)$ then $w\in C([0,T];H^1_0 (\ZOMq))\cap C^1([0,T];L^2(\ZOMq))$.

If $w_0=0$, $w_1=0$, $F=0$  and $f\in L^2(0,T;L^2(\Gamma))$ then $w\in C([0,T];L^2(\ZOMq))\cap C^1([0,T];H^{-1}(\ZOMq))$.

In every case $w$ depends continuously on the data in the specified spaces.
 
\end{Theorem}

These results justify the following definition of controllability:
\begin{Definition}\ZLA{eq:defiCONTROvisco}
Let $T>0$ and 
\[
\Lambda _{M}(T) f=w(T)\,,\qquad R_M(T)={\rm im}\, \Lambda _{M}(T)=\left \{ w(T)\,,\quad f\in L^2(0,T;L^2(\Gamma))\right \}\,.
\]
System~(\ref{eq:sistemaVERS4}) is controllable when the map $\Lambda _{M}(T)$ is surjective, i.e. when $R_M(T)=L^2(\ZOMq)$.

\end{Definition}

The result that we shall prove is:
\begin{Theorem}\ZLA{teo:CONTROLLprovatoDAVVERO}
Let the associated wave equation be controllable at time $2T$ and let $\ZEP>0$. Then system~(\ref{eq:sistemaVERS4}) is controllable at time $T+\ZEP$.

\end{Theorem}
 
 \begin{Remark}
 This observation will be important. The affine term $F$ in~(\ref{equA:FORMu;adiW}) does depend on $w_0$ and $w_1$ because of the integration by parts in the MacCamy trick. If the given equation is~(\ref{eq:dopoMacCamy}), with $F=0$ and $f=0$ then   formula~(\ref{equA:FORMu;adiW})  takes the form
 \begin{equation}
 \ZLA{eq:VolteRISOperAGGIU}
w(t)= R_+(t) w_0+\ZA^{-1}R_-(t)w_1+  
  \intt H(t-s) \left [R_+(s) w_0+\ZA^{-1}R_-(s)w_1 
 \right ]\ZD s 
 \end{equation} 
 \end{Remark}

\subsection{The direct inequality for Eq.~(\ref{eq:VolteRISOperAGGIU})}

Formula~(\ref{eq:DIREdisegPUREonde})  shows a ``hidden regularity'' of the wave equation, and this inequality  is called the ``direct inequality'' of the wave equation. We are going to prove an analogous resul for the solution of Eq.~(\ref{eq:VolteRISOperAGGIU}), i.e. we prove:
\begin{Theorem}\ZLA{Theo:DireINEAmemory} 
Let $T>0$. If $w_0\in H^1_0(\ZOMq)$ and $w_1\in L^2(\ZOMq)$ and let $ w $ solve 
Eq.~(\ref{eq:dopoMacCamy}) with $ f=0 $. Then $\zg_1w$ belongs to $ L^2(0,T;L^2(\Gamma))$ and  depends continuously on $w_0$, $w_1$, i.e. there exists $M$ such that
\begin{equation}\ZLA{eq:DireINEAmemory}
|\zg_1 w|^2_{L^2(0,T;L^2(\Gamma))}\leq M\left ( |w_0|^2 _{H^1_0(\ZOMq)} +|w_1|^2 _{L^2(\ZOMq)}+|F|^2_{L^1(0,T;L^2(\ZOMq))}\right )\,.
\end{equation}
\end{Theorem}

We give the proof in the case $F=0$ (the proof is easily adapted to $F\neq 0$).

The proof uses this property, that
\[
\mbox{if $\phi\in {\rm dom }\, A$ then $\zg_1\phi=-D^*A\phi$}\,.
\]

In order to prove Theorem~\ref{Theo:DireINEAmemory}
we introduce the notation $u(t)=R_+(t) w_0+\ZA^{-1}R_-(t) w_0$ and
\[
H_1(t)=\sum _{k=2}^{+\ZIN} \left ( \ZA^{-1}\right )^k L^{*k} =A^{-1}\sum _{k=2}^{+\ZIN} \left ( \ZA^{-1}\right )^{k -2}L^{*k} 
\]

 so that
\begin{equation}\ZLA{eq:RappreWsectIneDIre}
w(t)=u(t)+\ZA^{-1}\intt L(t-s) u(s)\ZD s+A^{-1}\intt H_1(t-s) u(s)\ZD s\,.
\end{equation}
Then,
\[
\zg_1 u(t)\,,\quad \zg_1\left (A^{-1}\intt H_1(t-s) u(s)\ZD s\right )=-D^*\left ( \intt H_1(t-s) u(s)\ZD s\right )
\]
are continuous functions of $w_0\in H^1_0(\ZOMq)$ and $w_1\in L^2(\ZOMq)$.

We study the first integral in(\ref{eq:RappreWsectIneDIre}),
\[
\ZA^{-1}\intt L(t-s) u(s)\ZD s=\ZA^{-1}\intt L(t-s)R_+(s) w_0\ZD s+A^{-1}\intt L(t-s) R_-(s) w_1\ZD s\,.
\]
The second term gives
\[
\zg_1\left (A^{-1}\intt L(t-s) R_-(s) w_1\ZD s\right )=-D^*\left ( \intt L(t-s) R_-(s) w_1\ZD s\right )\,,
\]
a continuous function of $w_1\in L^2(\ZOMq)$. We study the first integral. We recall that
\[
L(t)w=b R_-(t)w+\intt K(t-r)R_-(r)w\ZD r\,.
\]
We consider first
\begin{eqnarray*}
&&\ZA^{-1}\intt R_-(t-s) u(s)\ZD s=\\
&&=\ZA^{-1}\intt R_-(t-s) R_+(s) w_0\ZD s+A^{-1}\intt R_-(t-s) R_-(s) w_1\ZD s\;.
\end{eqnarray*}
The trace of the second addendum is treated as above. To handle the first addendum, we use
\[
R_-(\zt)R_+(r)=\frac{1}{2}\left ( R_-(r+\zt)-R_+(r-\zt)\right )
\]
so that
\[
\ZA^{-1}\intt R_-(t-s) R_+(s) w_0\ZD s=\frac{1}{2}t \left (\ZA^{-1}  R_-(t) w_0\right )+\frac{1}{2} \intt R_+(t-2s)\ZA^{-1}w_0\ZD s \,.
\]
The first addendum is the velocity term of the wave equation (even more regular, since $w_0\in H^1_0(\ZOMq)$) and the continuity of the trace follows from the properties of the wave equation. The same property holds also for $\zg_1\left (R_+(t-2s)\ZA^{-1}w_0\right )$ (say on the interval $(-T,T)$).

The convolution of these terms with $K$ retain the required properties.

\section{\ZLA{Sect:Controllability}The proof of controllability}
In this section we prove Theorem~\ref{teo:CONTROLLprovatoDAVVERO}. The proof is in two steps. In the first step  we prove that $R_M(T)$ is a closed subspace of $L^2(\ZOMq)\times H^{-1} (\ZOMq)$ and that  $R_M(T)^\perp$ is finite dimensional. In the second step we prove $R_M(T)^\perp=0$, hence controllability.
 
\subsection{The first step: $ R_M(T) $ is closed and  $R_M^\perp(T)$ is finite dimensional}

\begin{Theorem}
 \ZLA{teo:CodimFINITA} Let the associated wave equation be controllable at time $ 2T $. Then $ R_{ M}(T) $ is closed with finite codimension.
 \end{Theorem}
 \zProof 
 In the study of $ R_{M}(T) $ we use the notation
 \[ 
 u(t)=-\ZA\intt R_-(t-s)Df(s)\ZD s\,.
  \] 
  
  We fix any $ \zg<1/4 $. It is known that $ {\rm im}\, D\subseteq H^{1/2}(\ZOMq)\subseteq {\rm dom}(-A)^{\zg} $ and $ (-A)^{\zg} $ can be interchanged with  $ R_+(t) $ and $ R_-(t) $ and $ L(t) $. 
  
  We note that
 \begin{eqnarray*} 
&& \ZA^{-1}\intT L(t-s)u(s)\ZD s=-\intT L(t-s)\ints R_-(s-r)Df(r)\ZD r\,\ZD s=\\
 && (-A)^{-\zg}\intT L(t-s)\ints R_-(s-r)(-A)^\zg Df(r)\ZD r\,\ZD s\,.
  \end{eqnarray*}
  This is the composition of a continuous transformation with the compact transformation $ (A)^{-\zg} $. Hence it is a compact operator.
   For the same, and stronger, reasons  the map
   \[ 
f\ \mapsto K_T f=   \ \ZA^{-1}\intT L( T-s) u(s)\ZD s+A^{-1}\left [\sum _{k=2}^{+\ZIN} \left ( \ZA^{-1}\right )^{k-2} L^{*k}*u\right ] (T)
    \]
 is compact, from $ L^2(0,T;L^2(\ZOMq)) $ to $ L^2(\ZOMq) $.

 Then we have
 \[ 
 R _{M}(T) ={\rm im} \, \left ( \Lambda_0^1(T)+K_T\right )\;.
  \]
  The operator $ \Lambda_0^1(T) $ is surjective in $ L^2(\ZOMq)  $ by assumption while we proved that $ K_T $ is compact.
   
  Hence, $ R_M(T) $ is closed with finite codimension, as wanted.\zdia
  
  \subsection{The space $ R_M(T)^\perp $}
 
 We characterize $R_{M}(T)^\perp\subseteq L^2(\ZOMq) $:    
\[
\left (R_M(T)\right )^\perp =\left \{ \xi_0\in L^2(\ZOMq) \,,\quad \int_\ZOMq \xi_0(x) w(x,T)\ZD x =0\right \}\,.
\]
   This characterization will be applied also to the elements of $
  R_{M}(T)^\perp
  $ and we note that 
  \[
    R_{M}(T+\ZEP)^\perp\subseteq   R_{M}(T)^\perp\,.
  \]
   In this computation, closure of the reachable set has no interest, so that we can work with smooth controls. For example we can assume $f\in \mathcal{D}(\Gamma\times (0,T))$.

   We compute   $ \int_\ZOMq \xi_0(x)w(x,T)\ZD x$:
 
 \begin{align} 
 \nonumber &-\int_\ZOMq \xi_0  (x) \left [
  \ZA \intT R_-(T-s) D f(s)\ZD s+\right.\\
  \nonumber&+\left.  \ZA\intT H(T-s)\ints R_-(s-r)D f(r)\ZD r\,\ZD s
 \right ]\ZD x=\\
 \nonumber&=\left.- \intT\int_\Gamma f(r) D^*\ZA R_-(T-r)\xi_0  \ZD r\,\ZD\Gamma+\right.\\
   &\nonumber   +\intT\int_\Gamma f(r)D^*\ZA\int _0^{T-r}H(T-r-s)R_-(s)\xi_0  \ZD s\,\ZD \Gamma\,\ZD r=\\
  \ZLA{PriCOMPO}&=-\intT\int_\Gamma f(r) D^* A\left [ \ZA^{-1}\left (  R_-(T-r)\xi_0  + \int_0^{T-r}H(T-r-s) R_-(s)\xi_0  \ZD s\right )\right ]\ZD\Gamma\,\ZD r\,.
 \end{align}   

\begin{Remark}
Note that this is not a formal computation because the transformation $f\mapsto w$ is continuous.
\end{Remark}   
  
   If $ \xi_0   \perp R_M(T)$ then 
   
   \begin{equation}\ZLA{eq:TRAcciaNULLA}
 D^*A \left ( \ZA^{-1}R_-(r) \xi_0  + \intr  H(r-s) \ZA^{-1}R_-(s)  \xi_0   \ZD s \ZD\Gamma\,\ZD t \right)=0
   \end{equation} 
  Let
 \begin{eqnarray*}
&&\psi(t)= \ZA^{-1}R_-(t)\phi_1+ \intt  H(t-s)\ZA^{-1}R_-(s) \xi_0  \ZD s \,.
   \end{eqnarray*}
 We compare with~(\ref{eq:VolteRISOperAGGIU}) and we see that $\psi(t)$ solves 
   \begin{equation}\ZLA{equa:aggiuntaOLD}
\psi''=\Delta\psi+b\psi+\intt K(t-s)\psi(s)\ZD s\quad
\left\{\begin{array}{l}
\psi(0)=0\,,\ \psi'(0)= \xi_0  \,,\\
 \psi=0\ {\rm on}\  \partial\ZOMq   
 \end{array}\right.
   \end{equation}
   
   Note that $\xi_0\in L^2(\ZOMq)$ so that $\psi(t)\in C([0,T];H^1_0(\ZOMq))\cap C^1([0,T];L^2(\ZOMq))$.
   
   The interpretation of~(\ref{eq:TRAcciaNULLA}) is simple: when $ \xi_0 $ is ``smooth'', then $-D^*A=\zg_1$ and the direct inequality shows that $ \zg_1 $ is a continuous function of $ \xi_0\in L^2(\ZOMq) $, i.e.:
   \begin{Theorem}
   We have $ \xi_0\perp R_M(T) $ if and only if the solution of~(\ref{equa:aggiuntaOLD}) has the additional property
   \[ 
   \zg_1\psi(t)=0\ {\rm on}\ (0,T)\,.
    \]
   \end{Theorem}
   
   \subsection{The proof that $R_M(T+\ZEP)^\perp=0$}
   
   Let $\xi_0\perp R_M(T+\ZEP)^\perp$. We are going to prove $ \xi_0=0 $. We expand
   \begin{equation}\ZLA{eq:expansDIpsi}
\xi_0(x)=\ZSUno \phi_n(x)\xi_n\,,\qquad \{\xi_n\}\in l^2\,.   
   \end{equation}
 %%%%%%%%%%%%%%  
   The solution $\psi$ of system~(\ref{equa:aggiuntaOLD}) has the expansion
   \[
\psi(x,t)=\ZSUno \phi_n(x) \psi_n(t)\xi_n   
   \]
   where
   $\psi_n(t)$ solves
   \begin{equation}\ZLA{equaDIpsiN}
\psi_n''=-\zl_n^2\psi_n   +b\psi_n(t)+\intt K(t-s)\psi_n(s)\ZD s\,,\qquad \psi_n(0)=0\,,\quad \psi_n'(0)=1\,.
   \end{equation}
   The condition  $\xi_0\perp R_M(T+\ZEP)$    is the condition
   \begin{equation}\ZLA{condBorDOsvilu0}
\zg_1\psi(t)=
\ZSUno \left (\zg_1\phi_n\right )\xi_n \psi_n(t)   =0\,,\qquad 0<t<T+\ZEP\,.
   \end{equation}
   
   \begin{Remark}\ZLA{Remark:sullaCONVesre} This is a consequence of the direct inequality which  implies 
\[
\lim_N \sum _{n=1^N} \left (\zg_1\phi_n\right )\xi_n \psi_n(t)=\zg_1\psi
\]   
in $L^2(0,T+\ZEP;L^2(\Gamma))$.   \zdia
   \end{Remark}   
   %%%%%%%%%%%%%%%
    
   The goal is the proof that  equality~(\ref{condBorDOsvilu0}) implies $\xi_0=0$.

   In principle, it might be that the series in~(\ref{condBorDOsvilu0}) is a finite sum, i.e. that $\xi_n=0$ for large $n$.

We consider first the case that the series~(\ref{eq:expansDIpsi}) is a finite sum and then the case that it has infinitely many nonzero elements.

\subparagraph{The case $\xi_0=\sum _{n=1}^{N } \xi_n\phi_n $}
The sum cannot   have only one addendum, since otherwise we should have
\[
\zg_1\phi_{n_0}=0\ {\rm on}\ \Gamma
\] 
and $\phi _{n_0}$ is an eigenvector of $A$ and $\Gamma$ is the active part of $\partial\ZOMq$. It is known that this is not possible if there exists a time at which the wave equation is controllable. Even more, the terms with nonzero coefficients $\xi_n$ must belong to different eigenvalues, see~\cite{TAOcorrect,PandLIBRO}.

So, the sum must have at least two terms (which correspond to different eigenvalues) and we can assume $\xi_N\neq 0$. The fact that $\xi_0\perp R_M(T)$ implies
\begin{equation}\ZLA{equa:perpeIUNo}
\sum _{n=1}^{N } \xi_n\left (\zg_1 \phi_n\right )\psi_n(t)=0\,.
\end{equation}
Hence, also the second  derivative is zero and this, coupled with~(\ref{equa:perpeIUNo}), gives
\begin{equation}\ZLA{equa:perpeIDUE}
\sum _{n=1}^{N } \zl_n^2\xi_n\left (\zg_1 \phi_n\right )\psi_n(t)=0\,.
\end{equation}
We multiply~(\ref{equa:perpeIUNo}) with $\zl_N^2$ and we subtract from~(\ref{equa:perpeIDUE}). We get
\[
\sum _{n=1}^{N-1} \left (\zl_n^2-\zl_N^2\right )\xi_n\left (\zg_1 \phi_n\right )\psi_n(t)=0\,.
\]
If in this sum the nonzero coefficients $\left (\zl_n^2-\zl_N^2\right )\xi_n$ correspond to the same eigenvalue, this contradicts the previous observation. But, after a finite number of iteration of the procedure surely we obtain this case, which is not possible. Hence,  if $\xi_0\neq 0$ then the sum cannot be finite.

\subparagraph{Infinitely many nonzero entries}

The analysis of this case requires an intermediate step: we prove that $\xi_0$ is smoother then solely square integrable. In fact we prove:
\begin{Theorem}\ZLA{teo:regolaritaRTperp}
Let the wave equation be controllable at time $T$ and let $\ZEP>0$.
If $\xi_0\in L^2(\ZOMq)$ belongs to $R_M(T+\ZEP)^\perp$ then we have   
    $\xi_0\in {\rm dom}\, A$, i.e.
   \[
\xi_0(x)=\ZSUno \frac{\ZSI_n}{\zl_n^2}\phi_n(x)\,,\qquad \{\ZSI_n\}\in l^2\,.   
   \]
   \end{Theorem}
   
   We accept this theorem, whose proof is in the appendix, and we proceed to prove that $\xi_0=0$.

  We insert the special form of $\{\xi_n\}$ in~(\ref{condBorDOsvilu0}) and we find
  \[
  \ZSUno \left (\zg_1\phi_n\right )\frac{\ZSI_n}{\zl_n^2} \psi_n(t)   =0\,.
  \] 
   The observation in Remark~\ref{Remark:sullaCONVesre} implies that 
   \[
  \ZSUno \left (\zg_1\phi_n\right ) \ZSI_n \psi_n(t)  
  \] 
  is convergent. And so the following equality holds:
  \begin{eqnarray*}
  &&0=\frac{\ZD^2}{\ZD t^2}\ZSUno \left (\zg_1\phi_n\right )\xi_n \psi_n(t)=
  -\ZSUno \left (\zg_1\phi_n\right )\left (\zl_n^2\xi_n\right ) \psi_n(t)+\\
  &&+ \ZSUno \left (\zg_1\phi_n\right )
  \left [ b\psi_n(t)+\intt K(t-s)\psi_n(s)\ZD s\right]\xi_n=-\ZSUno \left (\zg_1\phi_n\right )\ZSI_n \psi_n(t)\;.
  \end{eqnarray*}
  
  This is the condition that 
  \[
\xi_1=\ZSUno \phi_n(x)\ZSI_n=  \ZSUno \phi_n(x)\left (\zl_n^2  \xi_n\right )\perp R_M(T)\,.
  \]
 So, using $\xi_0\perp R_M(T+\ZEP)$ we constructed a second element $\xi_1\perp R_M(T+\ZEP)$ and   the two elements $\xi_0$ and $\xi_1$ are linearly independent \emph{thanks to the fact that  (at least) two entries of $\xi_0$ which correspond to \emph{different eigenvalues} are nonzero.} 
 
 The new element 
 \[
\xi_1= \ZSUno \phi_n(x)\ZSI_n
 \]
 has the same properties as $\xi_0$ and so the procedure can be repeated. We get a third element $\xi_2\perp R_M(T+\ZEP)$,
 \[
\xi_2= \ZSUno \phi_n(x)\left (\zl_n^4  \xi_n\right )\in L^2(\ZOMq)
 \]
 and  the vectors $\xi_0$, $\xi_1$ and $\xi_2$ are linearly independent \emph{since   (at least) three entries of $\xi_0$ which correspond to \emph{different eigenvalues} are nonzero.}  
 
 The procedure can be iterated as many times  as we want, because we assumed that $\xi_0$ has infinitely many non zero entries (while every eigenvalue has finite multiplicity)  and we find that  ${\rm dim}\,R_M(T+\ZEP)^\perp=+\ZIN$. We proved already that this is false and so we get that any element  $
   \xi_0\perp R_M(T+\ZEP)$ has to be zero: $\xi_0=0$. This is the result that we wanted to achieve.
   
   \section{Appendix: the proof of Theorem~\ref{teo:regolaritaRTperp}}
    It is known that

   \[
{\rm dim}\, \ZOMq=d\ \implies m_0 n^{2/d}\leq \zl_n^2\leq M   n^{2/d}\,,\qquad m_0>0\,.
   \]
  In this proof we use the condition ${\rm dim}\, \ZOMq\leq 3$
which implies
  \begin{equation}\ZLA{eq:APPE:serieDAconvergere}
\left \{\zl_n^2\right \}\in l^2\ {\rm i.e.}\ \ZSUno \frac{1}{\zl_n^{4}}<+\ZIN
\end{equation}
   but it will be clear that this condition can be easily removed.   
   Furthermore we present the computation in the case $ b=0  $, only for simplicity of notations.    
   We shall see that this condition has no real effect on the computations.
   
   We use
   \begin{equation}\ZLA{eq:EquadiPsinAppe}
   \psi_n(t)=\frac{1}{\zl_n}\xi_n \sin\zl_n t + \intt \left [\frac{1}{\zl_n}\int_0^{t-s}K(r)\sin\zl_n(t-s-r)\ZD r\right ]\psi_n(s)\ZD s  \,.
    \end{equation}

We introduce the notations 
\[ 
S_n(t)=\sin\zl_n t\,,\qquad C_n(t)=\cos\zl_n t
 \]
 and $ L_n(t) $, the resolvent kernel of   the bracket in~(\ref{eq:EquadiPsinAppe}) (with the sign changed) so that
\begin{eqnarray}
\nonumber L_n &=& -\frac{1}{\zl_n}K*S_n+\frac{1}{\zl_n}(K*S_n)*L_n=\\
 \ZLA{eq:appe:SviluppddiLn}&=& -\frac{1}{\zl_n}
 K*S_n-
 \frac{1}{\zl_n^2}K^{*2} *S_n^{*2} 
 +\frac{1}{\zl_n^2}\left (K^{*2}  *S_n^{*^2}\right )*L_n\,.
\end{eqnarray}
 The first line of~(\ref{eq:appe:SviluppddiLn}) shows that 
 \begin{equation}\ZLA{eq:StomaLnT}
  |L_n(t)|\leq M/\zl_n  \ {\rm for} \ t\in (0,T) \,.
  \end{equation}
  
Due to the fact that the associated wave equation is controllable in time $ 2T $, hence also in larger times, we know that both $\{ \Psi_nS_n \} $ and $\{\Psi_n C_n \} $ where $\Psi_n=\zg_1\phi_n/\zl_n$ are Riesz sequences in $ L^2(0,T;L^2(\Gamma)) $ and in $ L^2(0,T+\ZEP;L^2(\Gamma)) $ and so the series
\[ 
\ZSUno \xi_n \Psi_nS_n\,,\qquad \ZSUno \xi_n\Psi_n C_n
 \]
   converge when $ \{\xi_n\}\in l^2 $.

Now we use
\[ 
\psi(x,t)=\ZSUno \phi_n(x)\psi_n(t)\xi_n \,,\qquad \psi_n(t)=\frac{1}{\zl_n} S_n(t)-\frac{1}{\zl_n} \left ( L_n*S_n\right )(t)\,.
 \]    
  So,  the condition of orthogonality to $ R_M(T) $ is
 \[
 \ZSUno\left \{ \xi_n \Psi_n S_n- \xi_n\Psi_n     \left (L_n*S_n\right )\right )=0\,. 
 \]
 This series converges and the equality holds  in $L^2(0,T+\ZEP;L^2(\Gamma))$ and, as we noted, the series $ \ZSUno \xi_n \Psi_n S_n$ converges too, so that we can write
 
 \[
 \ZSUno \xi_n \Psi_n S_n=\ZSUno\xi_n\Psi_n     \left (L_n*S_n\right )\,.
 \]
  
 We prove that  this function belongs to $H^1(0,T+\ZEP;L^2(\Gamma))$.   We formally compute termwise the derivative of the series on the right hand side and we prove that the resulting series converges in $L^2(0,T;L^2(\Gamma))$. In fact, the derivative is
   \begin{eqnarray} 
\nonumber&& \ZSUno \Psi_n\xi_n  \left (\zl_n  L_n*C_n\right )= -\ZSUno \Psi_n\xi_n  K*S_n*C_n-
  \\
 \nonumber && -\ZSUno\Psi_n \xi_n  \frac{1}{\zl_n }K^{*2}* S_n^{*2}*C_n+\\
\ZLA{eq:Appe:PAs1Se2} &&+\ZSUno \Psi_n \xi_n  \frac{1}{\zl_n }K^{*2}* S_n^{*^2}*C_n*L_n\,.
    \end{eqnarray}
 The first and second series on the right hand side converge  since
  \[ 
    S_n*C_n=\frac{1}{2}tS_n\,,\quad  S_n^{*2}*C_n=-\frac{1}{8} \left [t^2C_n(t)-\frac{1}{\zl_n}tS_n(t)\right ]\,.
     \]
The third series converges (even uniformly) since, using~(\ref{eq:StomaLnT}),  
     \begin{equation}  \ZLA{eq:SimadiLdivisoLAMBDA}
\left |\frac{1}{\zl_n}L_n\right |\leq \frac{M}{\zl_n^2}    \;. 
     \end{equation}
     Hence we have
     \[
     \ZSUno \xi_n \Psi_n S_n\in H^1(0,T+\ZEP;L^2(\Gamma)) \,.
     \]
    We combine with the fact that $\{\Psi_nS_n\}$, $\{\Psi_n C_n\}$ (and $\{\Psi_n e^{i\zl_n t}\}$) are Riesz sequences on the \emph{shorter} interval $(0,T)$ and we deduce (see~\cite[Chapt.~3]{PandLIBRO})
     \[
\xi_n=\frac{\ZDE_n}{\zl_n}\,,\qquad \{\ZDE_n\}\in l^2\,.     
     \]
     We replace this expression of $\xi_n$ and we equate the derivatives of both the sides. We get:
     \begin{eqnarray*}
      \ZSUno \ZDE_n \Psi_n C_n&=&\\
      %%%%%%%%%%
  &&     -\ZSUno \Psi_n\frac{\ZDE_n}{\zl_n}  K*S_n*C_n-
  \\
 \nonumber && -\ZSUno\Psi_n \frac{\ZDE_n}{\zl_n}  \frac{1}{\zl_n }K^{*2} *S_n^{*2}*C_n+\\
 &&+\ZSUno \Psi_n \frac{\ZDE_n}{\zl_n} \frac{1}{\zl_n }K^{*2}* S_n^{*^2}*C_n*L_n\,.
           \end{eqnarray*}
  Now  we see that the right hand side belong to $H^1(0,T;L^2(\Gamma))$. In fact, computing the derivatives termwise of the three series we get
  \begin{eqnarray}
&&  \ZLA{DeriSERI1} \ZSUno \Psi_n \ZDE_n   K* C_n^{*2}\,,\\
&&  \ZLA{DeriSERI2}           \ZSUno\Psi_n \ZDE_n   \frac{1}{\zl_n }K^{*2}  *C_n^{*2} * S_n\,,\\
 &&   \ZLA{DeriSERI3}           \ZSUno \Psi_n \ZDE_n  \frac{1}{\zl_n }K^{*2}* S_n *C_n^{*^2}*L_n\,.
  \end{eqnarray}   
 The series~(\ref{DeriSERI1}) and~(\ref{DeriSERI2})  converge  since
   \begin{eqnarray*}
   C_n^{*2}(t)&=&\frac{1}{2} \left (tC_n(t)+\frac{1}{\zl_n}S_n(t)\right )\,,\\
   S_n*C_n^{*2}&=&\frac{1}{8}\left [  \left ( t^2+\frac{1}{\zl_n^2}\right )S_n(t)-\frac{1}{\zl_n}t C_n(t)\right ]\,.
     \end{eqnarray*}
The series~(\ref{DeriSERI2}) and~(\ref{DeriSERI3}) converge, even uniformly, thanks to the inequality~(\ref{eq:SimadiLdivisoLAMBDA}).
     
  Hence we have
  \[
 \ZSUno \ZDE_n \Psi_n C_n\in H^1(0,T;L^2(\ZOMq))\quad \mbox{so that}\quad \ZDE_n=\frac{\ZSI_n}{\zl_n}  
  \]   
  hence
  \[
\xi_n=\frac{\ZSI_n}{\zl^2_n}  \,,
  \]   
   as we wanted to prove.

   \begin{Remark}
   The condition $ {\rm dim}\,\ZOMq\leq 3 $ has been used
when we replace  $ L_n(t) $ with its representation in the second line of~(\ref{eq:appe:SviluppddiLn}),   which has a coefficient $1/\zl_n^2$. Then we use   $ \{1/\zl_n^2\}\in l^2 $. If  $ {\rm dim}\,\ZOMq>3  $
   then we have $ \{1/(\zl_n^{2k}) \}\in l^2 $ provided $ k $ is sufficiently large. And we can get a factor $1/(\zl_n^{2k})  $ in~(\ref{eq:appe:SviluppddiLn}) by taking iterates of sufficiently high order. So, the condition ${\rm dim}\,\ZOMq\leq 3$ is easily removed.
   
   Also the condition $b=0$ it is easily removed: it is sufficient to replace $\zl_n$ with $\beta_n=\sqrt{\zl_n^2-b}$.\zdia
   \end{Remark}


\begin{thebibliography}{99}
 
  \bibitem{AvdoninPANDOLFIflux1}
S. Avdoni, L. Pandolfi:  
 {Simultaneous temperature and flux controllability for heat equations with memory.} 
  {Quarterly Appl. Math.,\/} {\bf 71}   339-368, 2013. 
  
    \bibitem{AvdoninPANDOLFIflux2}
    S. Avdonin,  L. Pandolfi:
  {Temperature and heat flux dependence/independence for  heat equations with memory.} 
    in``Time Delay Systems - Methods, Applications and New Trend'' (Ed.
R. Sipahi, T. Vyhlidal, P. Pepe, and S.-I. Niculescu)  Lecture Notes in Control and
Inform. Sci.~423, Springer-Verlag,   New York, 2012,  pp.~87-101. 


\ZBI{FuZHANG1}  X. Fu, J. Yong, X. Zhang:
 { Controllability and observability of the heat
equation with hyperbolic memory kernel.\/}
  {  J. Diff. Equations,\/}  {\bf 247}   2395-2439, 2009.
 
 \ZBI{TAOcorrect}A. Hassel, T. Tao:  Erratum for ``Upper and lower bounds for normal derivatives of Dirichlet eigenfunctions''. {\em Math. Res. Lett.\/} {\bf 17}   793-794, 2010.
 
 
  \ZBI{LeugeMOME} G.  Leugering:  On boundary controllability of viscoelastic systems. In {\em Control of partial differential equations\/} (Santiago de Compostela, 1987) Lecture Notes in Control and Inform. Sci., 114, Springer, Berlin   190-220, 1989.
  
  \ZBI{LIONSlibro} J.L. Lions:   \emph{Contr\^olabilit\`e exacte, perturbations et stabilization de syst\'emes distribu\`es.\/}
Vol.~1, Masson, Paris, 1988.

\ZBI{LoretiKOMORNIK} V. Komornik,  P. Loreti:  \emph{Fourier series in 
control theory.\/} Springer Monographs in Mathematics. Springer-Verlag, New York,   2005.
  
  \ZBI{Kim1993} J.U. Kim:  Control of a second-order integro-differential equation. {\em SIAM J. Control Optim. } {\bf 31}   101-110, 1993.
  
   \bibitem{PandViscoUnderTRACT}    
  L. Pandolfi:
 {Boundary controllability and source reconstruction in a viscoelastic string under external traction}.
  \emph{submitted for publication,  arXiv:1206.3034 J. Math. Analysis Appl.} (2013) DOI:10.1016/j.jmaa.2013.05.051
  
 
 
 \ZBI{PandAMO}  L. Pandolfi:
 {The controllability of the Gurtin-Pipkin equation: a cosine operator approach.} % {\em
 Appl.Math. and Optim.~52   143-165, 2005 (a correction is in { Appl. Math. Optim.\/}~64    467-468, 2011).

  \bibitem{PandDCDS2} %(MR2746431)
 L. Pandolfi,
 {  Riesz systems and an identification problem for heat equations with memory.\/} 
\emph{Discrete Continuous Dynam. Systems-S}   {\bf 4} (2011) 745--759.

\ZBI{PANDsurvey}L. Pandolfi: On-line input identification and application to Active Noise Cancellation. \emph{Annual Reviews in Control,}  vol. 34, 245-261, 2010. 

\ZBI{PandolfiTRACTION-deformation} L. Pandolfi: Traction, deformation and velocity of deformation in a viscoelastic string. Evolution equations Control Theory, {\bf 2}   471-493, 2013.

   \ZBI{PandolfiSHARP} L. Pandolfi: Sharp control time for viscoelastic bodies.  Submitted, arXiv:1305.1477. 
   
 \ZBI{PandLIBRO} L. Pandolfi: \emph{Distributed Systems with Persistent Memory: Control and Moment Problems.} Springer, in preparation.
 
 \ZBI{TuksnakWeiss} M. Tucsnak, G. Weiss:   \emph{Observation and control for operator semigroups. \/} Birkh\"auser, Basel, 2009.

 \end{thebibliography}
 \enddocument